\documentclass[12pt]{article}

\usepackage{latexsym}
\usepackage{amssymb}
\usepackage{verbatim}
\usepackage{amsthm}
\usepackage{amsmath}
\usepackage{yfonts}

\let\cal\mathcal

\let\cal=\mathcal      

\def\mcc{M\raise.5ex\hbox{c}C}
\def\mccarthy{M\raise.5ex\hbox{c}Carthy}

\def\eg{{\it e.g. }}
\def\ie{{\it i.e. }}

\def\h{{\cal H}}

\def\K{{\cal K}}
\def\M{{\cal M}}

\def\m{Mult}

\def\ga{\gamma}

\def\La{\Lambda}

\def\l{\lambda}

\def\vare{\varepsilon}

\let\i=\infty

\def\la{\langle}
\def\ra{\rangle}
\def\={\ = \ }    
\def\ot{\otimes}

\def\E{E_\l}

\def\kp{{k^\prime}}

\def\D{\mathbb D}

\def\be{\setcounter{equation}{\value{theorem}} \begin{equation}}
\def\ee{\end{equation} \addtocounter{theorem}{1}}
\def\beq{\begin{eqnarray*}}
\def\eeq{\end{eqnarray*}}

\def\vs{\vskip 5pt}

\def\bp{{\sc Proof: }}
\def\ep{{}{\hfill $\Box$} \vskip 5pt \par}

\def\bl{\begin{lemma}}
\def\el{\end{lemma}}
\def\bt{\begin{theorem}}
\def\et{\end{theorem}}
\def\bprop{\begin{prop}}
\def\eprop{\end{prop}}
\def\bd{\begin{definition}}
\def\ed{\end{definition}}
\def\br{\begin{remark}}
\def\er{\end{remark}}
\def\bexer{\begin{exercise}}
\def\eexer{\end{exercise}}
\def\bfig{\begin{figure}}
\def\efig{\end{figure}}

\numberwithin{equation}{section}

\usepackage[usenames,dvipsnames] {color}

\title{Pick Interpolation for free holomorphic functions}
\author{Jim Agler
\thanks{Partially supported by National Science Foundation Grant
DMS 1068830}
\\ U.C. San Diego\\ La Jolla, CA 92093
\and
John E. M\raise.5ex\hbox{c}Carthy
\thanks{Partially supported by National Science Foundation Grant DMS 
DMS 1300280}
\\ Washington University\\ St. Louis, MO 63130}

\def\be{\begin{equation}}
\def\ee{\end{equation}}

\def\ip#1#2{\langle #1,#2 \rangle}
\def\m{\mathbb{M}}
\def\md{\mathbb{M}^{[d]}}

\def\mn{\mathbb{M}_n}
\def\mnd{\mathbb{M}_n^d}
\def\mmd{\mathbb{M}_m^d}

\def\c{\mathbb{C}}
\def\D{\mathbb{D}}

\def\n{\mathbb{N}}
\def\cn{{\mathbb{C}^n}}
\def\k{\mathcal{K}}
\def\h{\mathcal{H}}

\def\L{\mathcal{L}}

\def\idcn{{\rm id}_{\mathbb{C}^n}}

\def\id#1{{\rm id}_#1}
\def\idd{{\rm id}}

\def\calv{\mathcal{V}}
 \def\calvccj{{\mathcal V}_{{\L(\c,\c^J)}}}

\def\calvhm{\mathcal{V}_{\L(\h,\M)}}

\def\vspace{\mathcal{V}}

\def\htwolo{{{\rm  H}_{L_1}^2}}
\def\htwolt{{{\rm  H}_{L_2}^2}}

\def\pd{\mathbb{P}^d}
\newcommand\gdel{G_\delta}

\def\si{s^{-1}}

\def\ot{\otimes}
\def\La{\Lambda}
\def\l{\lambda}

\def\Il{{\mathcal I}_\Lambda}
\def\Vl{V_\Lambda}
\def\Ml{{\mathbb M}_\Lambda}
\def\Mx{{\mathbb M}_x}

\renewcommand\={\ = \ }
\def\vs{\vskip 5pt}

\def\la{\langle}
\def\ra{\rangle}
\def\i{\infty}
\renewcommand\d{\delta}

\def\E{{\mathcal{E}}}

\def\bd{\begin{defin}}
\def\ed{\end{defin}}
\def\be{\begin{equation}}
\def\ee{\end{equation}}
\def\begm{\left(\begin{matrix}}
\def\endm{\end{matrix}\right)}

\def\E2{E^{[2]}}
\def\tr{{\rm tr}}
\def\Np{N}
\def\Nz{{N_0}}

\def\cno{{\mathcal N}_1}
\def\cnt{{\mathcal N}_2}
\def\cmo{{\mathcal M}_1}
\def\cmt{{\mathcal M}_2}

\def\clo{{\mathcal L}_0}
\def\higd{{H^\i(\gdel)}}
\def\bhigd{{\rm ball}(\higd)}
\def\kp{{\k}^\prime}
\def\cj{{\c^J}}
\def\cm{{\c^m}}
\def\op{\oplus}

\def\gv{\textswab{V}}
\def\ga{\textswab{A}}
\def\varga{{\rm Var}(\ga)}

\def\hinfgdel{H^\infty(\gdel)}
\def\hinfo{H_1^\infty(\gdel)}
\def\hinfu{H^\infty(U)}
\def\hinfuo{H_1^\infty(U)}

\theoremstyle{definition}
\newtheorem{theorem}[equation]{Theorem}
\newtheorem{defin}[equation]{Definition}
\newtheorem{lem}[equation]{Lemma}
\newtheorem{lemma}[equation]{Lemma}
\newtheorem{prop}[equation]{Proposition}

\newtheorem{remark}[equation]{Remark}

\def\bl{\begin{lem}}
\def\el{\end{lem}}
\def\bt{\begin{theorem}}
\def\et{\end{theorem}}

\begin{document}

\bibliographystyle{plain}

\maketitle

\section{Introduction}
\label{seca}

A holomorphic function in $d$ (scalar) variables behaves  locally like a  polynomial.
Given such a function $\phi$, one can evaluate it also on $d$-tuples of commuting matrices
whose joint spectrum lies in the domain of $\phi$.
A free holomorphic function  is a generalization of this notion, where the matrices are no longer required to commute. 
We replace polynomials by free polynomials, \ie polynomials in non-commuting variables,
and consider functions that are locally limits of free polynomials.

To make this precise, let us first define
$\mnd$ to be the set of all $d$-tuples of $n$-by-$n$ complex matrices, and $\md = \cup_{n=1}^\i \mnd$. Let $\pd$ denote the algebra of free polynomials in $d$ variables.
A {\em graded function} defined on a subset of $\md$ is a function $\phi$ with the property that
if $x \in \mnd$, then $\phi(x) \in \mn$.
\begin{defin}\label{defa3}
An {\em nc-function} is a graded function $\phi$ defined on a
set $D \subseteq \md$ such that

i) If $x,y, x\oplus y \in D$, then $\phi(x \oplus y) = \phi(x) \oplus \phi(y)$.

ii) If $s \in \mn$ is invertible and $x, s^{-1} x s \in D \cap \mnd$ then $\phi(\si x s) = \si \phi(x) s$.

\end{defin}

Nc-functions have been studied for a variety of reasons: by Taylor \cite{tay73}, in the context of the functional calculus 
for non-commuting operators;
 Voiculescu \cite{voi4,voi10}, in the context of free probability;
Popescu \cite{po06,po08,po10,po11}, in the context of extending classical function theory to 
$d$-tuples of bounded operators;
 Ball, Groenewald and Malakorn \cite{bgm06}, in the context of extending realization formulas
from functions of commuting operators to functions of non-commuting operators;
Alpay and Kalyuzhnyi-Verbovetzkii \cite{akv06} in the context of
realization formulas for rational functions that are $J$-unitary on the boundary of the domain;
Helton \cite{helt02} in proving  positive matrix-valued functions are sums of squares;
 and
Helton, Klep and McCullough \cite{hkm11a,hkm11b}
and Helton and McCullough \cite{hm12} in the context of developing a descriptive theory of the domains on
which LMI and semi-definite programming apply.
Recently, Kaliuzhnyi-Verbovetskyi  and Vinnikov have written a monograph
\cite{kvv12}
that gives a panoramic view of the developments in the field to date, and establishes their
Taylor-Taylor formula for nc-functions.

There are two topologies that we wish to consider on $\md$. The first is called the {\em  disjoint union topology}: a set $U$ is open in the disjoint union topology if and only if $ U \cap \mnd$ is open for every $n$. This topology is too fine for some purposes; for example, compact sets must have a bound on the 
size of the matrices that they contain. The other topology we wish to consider is the {\em free topology},
which is most conveniently defined by giving a basis.
A basic free open set in $\md$ is a set of the form
\[
\gdel \= 
\{ x \in \md : \| \d(x) \| < 1 \},
\]
where $\delta$ is a $J$-by-$J$ matrix with entries in $\pd$.
We define the free topology to be the topology on $\md$ which has as a basis all the sets $\gdel$,
as $J$ ranges over the positive integers, and the entries of $\delta$ range over all polynomials in $\pd$.
(Notice that $ G_{\delta_1} \cap G_{\delta_2} = G_{\delta_1 \oplus \delta_2}$, so these sets do form the basis of a topology). The free topology is a natural topology  when considering semi-algebraic sets.
\begin{defin}
\label{defa4}
A free holomorphic function $\phi$ on a free open set $U$ in $\md$  is a function that is locally  
a bounded nc-function, \ie for every $x$ in $U$ there is a basic free open set $\gdel \subseteq U$ that
contains $x$ and such that $\phi |_{\gdel}$ is a bounded nc-function.
\end{defin}

It is a principal result of \cite{amfree} that free holomorphic functions are locally approximable by
free polynomials (see Theorem~\ref{thmb1} below).

The main result of this note is a criterion for solving a Pick interpolation problem
on a basic open set, Theorem~\ref{thmpick} below, and its  generalization to extending
bounded
free holomorphic functions off free varieties, Theorem~\ref{thme1}.

Let $\hinfu$ denote the bounded free holomorphic functions on a free open set $U$ with the supremum norm, and let $\hinfuo$ denote the closed unit ball 
 of $\hinfu$.
For $ 1 \leq i \leq N$, let $\l_i \in \gdel \cap \m^d_{n_i}$ and let $w_i \in \m_{n_i}$.
The Pick problem is to determine whether or not there is a function in $\hinfuo$ that maps
each $\l_i$ to the corresponding $w_i$.

Note first that if $U$ is closed under direct sums, then 
by letting $\La = \oplus_{i=1}^N \l_i$ and $W = \oplus_{i=1}^N w_i$, the original $N$ point problem
 is the
same as solving the one point Pick
problem of mapping $\La$ to $W$. 
Secondly, unlike in the scalar case, one cannot always solve the Pick problem if one drops the norm constraint. 
For example, no holomorphic function maps
$$\begin{pmatrix} 0 & 1 \\
0 & 0
\end{pmatrix}
\qquad
{\rm to} \qquad
\begin{pmatrix} 1 & 0 \\
0 & 0
\end{pmatrix}.
$$

To state the theorem, let us make the following definitions for  $\La$ in $\md$.
Define
\[
\Il \= \{ p \in \pd\, : \, p(\La) = 0 \}
\]
and
\[
\Vl \= \{ x \in \md \, : \, p(x) = 0 \ {\rm whenever\ } p \in \Il \} .
\]
Let 
\[
\Ml \= \{ p(\La) \, : \, p \in \pd \}.
\]
Note that since $\Ml$ is a finite dimensional vector space, it is closed.
\bt
\label{thmpick}
Let $\La \in \gdel\cap \mnd$ and $W \in \m_n$.
There exists a function $\phi$ in the closed unit ball of $H^\i(\gdel)$ such that 
$\phi(\La) = W$
if and only if 

(i) $W \in \Ml$, so there exists $p_0 \in \pd$ such that $p_0(\La) = W$.

(ii) $\sup \{ \| p_0(x) \| : x \in \Vl \cap \gdel \} \ \leq \ 1$.
\et
We prove this theorem in Section~\ref{secc}. Note that when $d=1$, the question
of whether $p_0$ can be found satisfying $p_0(\La) = W$ can be resolved by
looking at the Jordan canonical form of $\La$. In this basis, the algebra $\Ml$ has a straightforward description. When $d > 1$, the determination of $\Ml$ is more delicate; generically\footnote{For example,
 if $\La^1$ has $n$ distinct eigenvalues and $\La^2$ has no non-zero entry when the matrix
is expressed in the basis
given by the eigenvectors of $\La^1$.},
 however,
the algebra $\Ml$ will be all of $\m_n$.

In Section~\ref{secd}, we give a description in Theorem~\ref{thmd1}
of all the solutions of a (solvable) Pick problem  ---
this is called the Nevanlinna problem. Our approach is indebtted to  the solution in the scalar case by
J.~Ball, T.~Trent and V.~Vinnikov \cite{btv01}.

Theorem~\ref{thmpick} has a remarkable corollary. Suppose $\ga$
 is an algebra in $\pd$, and let $\gv = \varga$ be given by

\be
\label{eqe1}
\gv \= \{ x \in \m^d : p(x) = 0 \ \forall \ p \in \ga \}.
\ee
If $\Lambda$ is in $\gv$, then $\ga \subseteq \Il$, and
$\Vl \subseteq \gv$.
Let $U$ be a free open set in $\md$;
we shall say that a function $f$ defined on $ \gv \cap U$ is free holomorphic if, for every point
$x$ in $ \gv \cap U$ there is a basic free open set $\gdel \subseteq U$  containing $x$ and a free holomorphic function
$\psi$ 
defined on $\gdel$ such that $\psi |_{\gv \cap \gdel} = f  |_{\gv \cap \gdel}$.

In the scalar case, every holomorphic function defined on an analytic variety inside a domain of holomorphy extends to a holomorphic function on the whole domain, by a celebrated theorem of H.~Cartan \cite{car51}. The geometric conditions that guarantee that  all bounded holomorphic functions extend to be bounded on the whole domain have been investigated by Henkin and Polyakov \cite{henpol84} and Knese \cite{kn10b}; however, 
even when bounded extensions exist, the extension is almost never isometric \cite{agmc_vn}.
But in the matrix case,
any bounded free holomorphic function  on $\gv \cap \gdel$ does extend to a free holomorphic function on
$\gdel$  with the same norm.

\bt
\label{thme1}
Let $\gv$ be as in \eqref{eqe1}, and let
 $\delta$ be a matrix of free polynomials  such that
$\gv \cap \gdel$ is non-empty. Let $f$ be a bounded  free holomorphic function 
defined on $\gdel \cap \gv$.
Then there is a free holomorphic function $\phi$ on $\gdel$
that extends $f$ and such that
\be
\label{eqe2}
\| \phi \|_{H^\i (\gdel)} \= \sup_{x \in \gv \cap \gdel} \| f(x) \|
\ee
\et
We prove this in Section~\ref{sece}. In Section~\ref{secf} we give some applications.

The definition 
\eqref{eqe1} naturally leads one to ask what the ideal of $\gv$, the set
\[
I_\gv \= \{ p \in \pd : p(x) = 0 \ \forall \ x \in \gv \} ,
\] 
is.
In the complex case, the answer is simpler than in the scalar case, at least
if $\ga$ is finitely generated. In \cite{helmccu04}, Bergman, Helton and McCullough proved 
that $I_\gv$ is the smallest ideal containing $\ga$, provided this ideal is finitely generated.
The real (self-adjoint) case is more subtle --- see \eg \cite{chmn13}.

\section{Background material}
\label{secnb}

We shall need some results from \cite{amfree}.  The first we have already referenced:
\bt
\label{thmb1}
Let $D$ be a free domain and let $\phi$ be a graded function defined on $D$. Then $\phi$ is a free holomorphic function if and only if $\phi$ is locally approximable by polynomials.
\et

The second, \cite[Thm 8.1]{amfree}, says that a function is in $\hinfo$ if and only if it has a free $\d$-realization.

\begin{defin}
\label{defnb1}
Let $\phi$ be a graded function on $\gdel$, where $\delta$ is a $J$-by-$J$ matrix of free
polynomials.
A free $\d$-realization of $\phi$ is a Hilbert space $\L$, an isometry $V : \c \oplus (\c^J \otimes\L )\to \c \oplus (\c^J \otimes\L)$ that 
can be written
\[
V \= 
\bordermatrix{
&\c & \c^J \otimes\L \cr
\c& A &  B\cr
 \c^J \otimes\L&C & D},
\]
and such that
\beq
\lefteqn{
\phi(x) \= \id{\cn} \ot A \ + } \\
&
 (\idcn \otimes B) (\d(x) \otimes \id{{\L}} )[ \idcn \otimes \id{{\c^J \otimes \L}} - (\idcn \otimes D) (\d(x) \otimes \id{{\L}} )]^{-1}
\idcn \otimes C
\eeq
for all $x \in \gdel \cap \mnd$.

We call $\phi$ the transfer function of $V$ (where $\d$ is understood).
\end{defin}
\bt
\label{thmb2}
Let $\phi$ be a graded function on $\gdel$. Then $\phi$ is in $\hinfo$ if and only if $\phi$
has a free $\d$-realization.
\et

The third is a Montel theorem.
\bt
\label{thmb3}
Let $(\phi_i)_{i=1}^\i$ be a sequence in $\hinfuo$. Then there is a subsequence 
$(\phi_{i_j} )_{j=1}^\i$ and a function $\phi$ in $\hinfuo$ such that
$(\phi_{i_j} )_{j=1}^\i$ converges to $\phi$ locally uniformly on $U$ in the disjoint
union topology.
\et

\section{Proof of Theorem~\ref{thmpick}}
\label{secc}

Let $E = \Vl \cap \gdel$, and let 
\[
E^{[2]} \=
\{ (x,y) \, : \, x,y \in \Vl \cap \gdel \cap \mmd,\ {\rm for\ some\ } m \}.
\]


Let us start with some lemmata.

\bl
\label{lemb1}
Let $\La, x  \in \md$. The following are equivalent: 

(i) $ x \in \Vl$.

(ii) There is a homomorphism $\alpha : \Ml \to \Mx$ such that $\alpha(\La^r) = x^r$ for  $ r = 1, \dots, d$.

(iii) The map $p(\La) \mapsto p(x)$ is a well-defined map from $\Ml$ to $\Mx$.

(iv) The map $p(\La) \mapsto p(x)$ is a completely bounded homomorphism. 
\el
\bp
The equivalence of {\em (i) - (iii)} is by definition. That {\em (iii)} is equivalent to {\em (iv)} 
is because every bounded homomorphism defined on a finite dimensional space is automatically
completely bounded \cite{pau02}.
\ep
\bl
\label{lemc1}
Let $\phi$ be in  $H^\i(\gdel)$. Then there exists a polynomial $p_0 \in \pd$ so that
\be
\phi(x) \= p_0(x) \qquad \forall x \in  \Vl \cap \gdel .
\label{eqc1}
\ee
\el
\bp
By Theorem~\ref{thmb1}, the free function $\phi$ can be uniformly approximated on a free
neighborhood of $\La$ by free polynomials. In particular, since $\Ml$ is closed, there is a polynomial
$p_0$ such that $\phi(\La) = p_0(\La)$.

Fix $x \in \Vl \cap \gdel$. By another application of the same theorem, there is a free polynomial
$p_1$ such that $\phi(\La \oplus x) = p_1(\La \oplus x)$.
Therefore $p_0(\La) = p_1(\La)$, so by the definition of $\Vl$, we also have
$p_0(x) = p_1(x)$. Therefore \eqref{eqc1} holds, as desired.
\ep
We let  $\calv$  denote the vector space
of nc-polynomials on $E$, where we identify polynomials that agree on $E$;
and   we let $\calvhm$  denote the vector space
of  $\L(\h,\M)$-valued  nc-polynomials on $E$.
As any such polynomial on $E$ is uniquely determined by its values on $\La$,
the space of such functions is finite dimensional, if
 $\h$ and $\M$  are finite dimensional.

Consider the following vector spaces of functions on $\E2$, where all sums are over a
finite set of indices:

\beq
H(E) &\= &
\{ h(y,x) = \sum g_i(y)^* f_i(x) : f_i, g_i \in \pd \}
\\
R(E) &=& \{ h \in H(E) : h(x,y) = h(y,x)^* \} \\
C(E) &=& \{ h(y,x) = \sum u_i(y)^* [ \idd - \delta(y)^* \delta(x) ] u_i(x) \, : \\
&& \qquad u_i \ {\rm
is\ } \L(\c,\c^J)-{\rm valued\ nc\ polynomial} \} \\
P(E) &=& \{ h(y,x) = \sum f_i(y)^* f_i(x) : f_i \in \pd \}
\eeq
We topologize $H(E)$ with the norm 
\[
\| h (y,x) \| \= \| h(\La, \La ) \| .
\]

\begin{lemma}\label{lemc2}
Let $\h, \M$ be  finite dimensional Hilbert spaces, and 
let $F(y,x)$ be an arbitrary graded $\L(\M)$-valued function on $E^{[2]}$.
Let $\Nz
= \dim(\calvhm)$.
Then
if $G$ can be represented in the form
\[
G(y,x) = \sum_{i=1}^m g_i(y)^*F(y,x)g_i(x),\qquad\ (x,y) \in \E2
\]
where $m \in \n$ and $g_i \in \calvhm$ for $i=1,\ldots,m$, then $G$ can be represented in the form
\be\label{5.12}
G(y,x) = \sum_{i=1}^\Nz f_i(y)^*F(y,x)f_i(x),\qquad\ (x,y) \in \E2
\ee
where $f_i \in \calvhm$ for $i=1,\ldots, \Nz$.
\end{lemma}
\bp
Let $\langle e_l(x)\rangle_{l=1}^\Nz$ be a basis of $\calvhm$. For each $i = 1,\ldots,m$, let
\[
g_i(x) = \sum_{l=1}^\Nz c_{il}\ e_l(x).
\]
Form the $m \times \Nz $ matrix $C=[c_{il}]$. As $C^*C$ is an $\Nz \times \Nz$ positive semidefinite matrix, there exists an $\Nz \times \Nz$ matrix $A = [a_{kl}]$ such that $C^*C = A^*A$. This leads to the formula,
\[
\sum_{i=1}^m \overline{c}_{il_1}c_{il_2} = \sum_{k=1}^\Nz \overline{a}_{kl_1}a_{kl_2},
\]
valid for all $l_1,l_2 = 1,\ldots, \Nz$. If  $ (x,y) \in  \E2$, then
\beq
G(y,x) &\=& \sum_{i=1}^m g_i(y)^*F(y,x)g_i(x)\\
&=&\sum_{i=1}^m \big(\sum_{l=1}^\Nz c_{il}\ e_l(y)\big)^*F(y,x)\big(\sum_{l=1}^\Nz c_{il}\ e_l(x)\big)\\
&=&\sum_{l_1,l_2=1}^\Nz (\sum_{i=1}^m \overline{c}_{il_1}c_{il_2})e_{l_1}(y)^*F(y,x)e_{l_2}(x)\\
&=&\sum_{l_1,l_2=1}^\Nz (\sum_{k=1}^\Nz \overline{a}_{kl_1}a_{kl_2})e_{l_1}(y)^*F(y,x)e_{l_2}(x)\\
&=&\sum_{k=1}^\Nz \big(\sum_{l=1}^\Nz a_{kl}\ e_l(y)\big)^*F(y,x)\big(\sum_{l=1}^\Nz a_{kl}\ e_l(x)\big).
\eeq
This proves that \eqref{5.12} holds with $f_i = \sum_{l=1}^\Nz a_{il}\ e_l$.
\ep
 
\bl
\label{lemc3}
$C(E)$ is closed.
\el
\bp
By Lemma~\ref{lemc2}, every element in $C(E)$ can be represented in the form
\be
\label{eqc2}
 \sum_{i=1}^\Nz u_i(y)^* [ \idd - \delta(y)^* \delta(x) ] u_i(x) ,
\ee
where $\Nz = \dim {\cal V}_{\L(\c,\c^J)}$.
Suppose a sequence of elements of the form \eqref{eqc2}
approaches some $h$ in $H(E)$ at the point $(\La,\La)$:
\[
\sum_{i=1}^\Nz u^{(k)}_i(\La)^* [ \idd - \delta(\La)^* \delta(\La) ] u^{(k)}_i(\La) 
\to
h(\La,\La)\ {\rm as\ } k \to \i .
\]
Since $\La \in \gdel$, there is a constant $M$ such that, for each $i$ and $k$,
\[
\|  u^{(k)}_i(\La)  \| \leq M .
\]
Passing to a subsequence, one can assume that each $ u^{(k)}_i(\La) $ converges to some
$u_i(\La)$ (since $u_i^{(k)}$ is a graded $\L(\c,\c^J)$ valued function and $J < \i$).
By Lemma~\ref{lemb1}, we have 
\[
u^{(k)}_i(x) \to 
u_i(x)  \qquad \forall x \in E .
\]
Therefore, for all $(x,y) \in \E2$, we have
\beq
\sum_{i=1}^\Nz u^{(k)}_i(y)^* [ \idd - \delta(y)^* \delta(x) ] u^{(k)}_i(x) 
&\to&
\sum_{i=1}^\Nz u_i(y)^* [ \idd - \delta(y)^* \delta(x) ] u_i(x) 
\\
&=& h(y,x) .
\eeq
\ep

\bl
\label{lemc4}
We have $P(E) \subseteq C(E)$.
\el
\bp
We have
\be
\label{eqc3}
f(y)^*f(x) - \sum_{k=0}^{m-1} f(y)^* \delta(y)^{k* }[ \idd - \d(y) \d(x) ]
\d(x)^k f(x) 
\=
f(y)^* \d(y)^m \d(x)^m f(x) .
\ee
As $m \to \i$, the right-hand side of \eqref{eqc3} goes to zero for every $(x,y) \in \E2$.
Since $C(E)$ is closed by Lemma~\ref{lemc3}, this proves that $f(y)^* f(x) \in C(E)$, and
hence so are finite sums of this form.
\ep

\bl
\label{lemc5}
Suppose  
$\sup \{ \| p_0(x) \| : x \in E \} \ \leq \ 1$.
Then the  function $$
h(y,x) \= \idd - p_0(y)^* p_0(x)
$$ is in $C(E)$.
\el
\bp
This will follow from the Hahn-Banach theorem \cite[Thm. 3.3.4]{rud91}  if we can can show that
$L(h(y,x)) \geq 0$ whenever 
\be
\label{eqc4}
L \in R(E)^* \qquad {\rm and }\qquad L(g) \geq 0 \ \forall\, g \in C(E) .
\ee 

Assume \eqref{eqc4} holds, and 
define $L^\sharp \in H(E)^*$ by the formula
\[
L^\sharp (h(y,x)) = L(\frac{h(y,x) + h(x,y)^*}{2})+iL(\frac{h(y,x) - h(x,y)^*}{2i}),
\]
and then define  sesquilinear forms  on $\calv$ and $\calvccj$ by the formulas
\begin{eqnarray*}
\langle f,g\rangle_{L_1} &\=& L^\sharp(g(y)^*f(x)),\qquad f,g \in \calv \\
\langle F,G\rangle_{L_2} &\=& L^\sharp(G(y)^*F(x)),\qquad F,G \in \calvccj.
\end{eqnarray*}

Observe that Lemma \ref{lemc4} implies that $f(y)^*f(x) \in C(E)$ whenever $f \in \calv$ or  $\calvccj$. Hence, \eqref{eqc4} implies that
$
\ip{f}{f}_{L_1} \ge 0
$
for all $f \in \calv$, and
$
\ip{F}{F}_{L _2}\ge 0
$
for all $F \in \calvccj$,
 i.e., $\ip{\cdot}{\cdot}_{L_1}$ and  $\ip{\cdot}{\cdot}_{L_2}$ are pre-inner products on $\calv$ and
$\calvccj$ respectively. 


To make them into inner products, choose $\vare > 0$ and define
\begin{eqnarray}
\label{eqc5}
\langle f,g\rangle_{1} &\= L^\sharp(g(y)^*f(x)) + \vare\,  \tr (g(\La)^* f(\La) ),\quad &f,g \in \calv \\
\label{eqc6}
\langle F,G\rangle_{2} &\= L^\sharp(G(y)^*F(x)) +  \vare\,  \tr (G(\La)^* F(\La) ),\quad &F,G \in \calvccj.
\end{eqnarray}

We let $\htwolo$ and  $\htwolt$ denote the Hilbert spaces
 $\calv$ and
 $\calvccj$ equipped with the  inner products \eqref{eqc5} and \eqref{eqc6}.

The fact that $L$ is non-negative on $C(E)$ means that
\be
\label{eqc7}
\langle F,F\rangle_{L_2} \ \geq \ 
\langle \d F,\d F\rangle_{L_2}
\ee
for all $F$ in $\calvccj$.
Since $\| \d(\La) \| < 1$, we also have
\be
\label{eqc8}
\tr ( F(\La)^* F(\La) ) \ > \ 
\tr ( F(\La)^* \d(\La)^* \d(\La) F(\La) ),
\ee
if $F \neq 0$, 
and combining \eqref{eqc7} and \eqref{eqc8}
we get that 
 multiplication by $\delta$ is a strict contraction
on $\htwolt$. 

Let $M$ denote the $d$-tuple of multiplication by the coordinate functions $x^r$ on $\htwolo$.
We have just shown that $\| \d(M) \| < 1$, so $M$ is in $\gdel$. As $M$ is also in $\Vl$,
we have that $M$ is in $E$.  Therefore $\| p_0(M) \| \leq 1$, by hypothesis.
Therefore
\[
\idd - p_0(M)^* p_0(M) \ \geq \ 0 ,
\]
and so for all $f$ in $\calv$ we have
\beq
\lefteqn{
L^\sharp(f(y)^*f(x)) + \vare\,  \tr (f(\La)^* f(\La) )
\ \geq \ } \\ &&
L^\sharp(f(y)^* p_0(y)^* p_0(x) f(x)) + \vare\,  \tr (f(\La)^* p_0(\La)^* p_0(\La) f(\La) ) .
\eeq
Letting $f$ be the function $1$ and letting $\vare \to 0$, we get
\[
L(\idd - p_0(y)^* p_0(x) ) \ \geq \ 0,
\]
as desired.
\ep

We can now prove the theorem.

{\sc Proof of Theorem~\ref{thmpick}:}
(Necessity). Condition {\em (i)} follows from Lemma~\ref{lemc1}. Condition {\em (ii)} follows
because $p_0(x) = f(x)$ for $x \in \Vl \cap \gdel$, and $f$ is in the unit ball of $H^\i(\gdel)$,
so $\| f(x) \| \leq 1$ for every $x$ in $\gdel$.

\vs

(Sufficiency). Suppose {\em (i)} and {\em (ii)} hold.
By Lemma~\ref{lemc5},  the function $$
h(y,x) \= \idd - p_0(y)^* p_0(x)
$$ is in $C(E)$.
By Lemma~\ref{lemc2}, there is some positive integer $\Np \leq {\rm dim}(\calv_{\L(\c^n,\c^n)})$ and an
$\L(\c,\c^{J\Np})$-valued nc polynomial $u$ such that, for $x,y \in E \cap \m_n^d$,
\begin{eqnarray}
\nonumber
h(y,x) &\=& \idd_{\c^n} - p_0(y)^* p_0(x) \\
&\=&
u(y)^* \left[ \idd_{\c^{nJ\Np}} - (\d(y)^* \otimes \idd_{\c^{\Np}}) \  (\d(x) \ot \idd_{\c^{\Np}} ) \right] u(x) .
\label{eqc10}
\end{eqnarray}
Replace $x$ in \eqref{eqc10} with $s x s^{-1}$ where $s$ is invertible in $\mn$ and
$sx s^{-1}$ is in $\gdel$ to get
\begin{eqnarray}
\label{eqc11}
\lefteqn{s - p_0(y)^* s p_0(x) \=  }\\
&&u(y)^* \left[s \ot \idd_{\c^{J\Np}} - (\d(y)^* \otimes \idd_{\c^{\Np}}) \ s\ot \idd_{\c^{J\Np}} (\d(x) \ot \idd_{\c^{\Np}}) \right] u(x) .
\nonumber
\end{eqnarray}
Equation~\eqref{eqc11} is true for all $s$ in a neighborhood of the identity, and as linear combinations
of such elements span $\m_n$, we get
that \eqref{eqc11}  actually holds for all $s$ in $\mn$.
For $k=1,\ldots,n$, define $\pi_k:\c^n \to \c$ by the formula
\[
\pi_k(v) = v_k,\qquad v = (v_1,\ldots,v_n) \in \c^n.
\]
Letting $s=\pi_l^*\pi_k$ in \eqref{eqc11} and applying to $v$ and taking the inner product with $w$,
where $v$ and $w$ are in $\c^n$, leads to
\begin{eqnarray}
\label{eqc12}
\lefteqn{
\la [ \pi_l^*\pi_k - p_0(y)^* \pi_l^*\pi_k p_0(x) ] v, w \ra \=
}\\
&&\la [ \pi_l^*\pi_k \ot \idd - (\d(y)^* \ot \idd )(\pi_l^*\pi_k \ot \idd)(\d(x) \ot \idd  )]
u(x) v, u(y) w \ra .
\nonumber
\end{eqnarray}

For each $v \in \c^n$ 
define vectors $p_v$  and $q_v$ in $\c^{n(1+\Np J)}$ by
\beq
p_v &\=&
\begin{bmatrix}
\id{\cn} \\
[\d(\La) \ot \id{{\c^{\Np}}}  ]u(\La) 
\end{bmatrix}
\ v 
\\
q_v &\=&
\begin{bmatrix}
p_0(\La) \\
u(\La) 
\end{bmatrix}
\ v .
\eeq
For each $1 \leq k \leq n$, define vectors $p_{k,v}$ and $q_{k,v}$ in $\c^{1 + \Np J}$ by
\beq
p_{k,v} &\=&  [\pi_k \ot \idd_{\c^{1 + \Np J}} ] p_v\\
q_{k,v} &\=&  [\pi_k \ot \idd_{\c^{1 + \Np J}} ] q_v.
\eeq
Then \eqref{eqc12}, with $\La$ in place of both $x$ and $y$, becomes
\be
\label{eqc31}
\la p_{k,v} , p_{l,w} \ra \=
\la q_{k,v} , q_{l,w} \ra \qquad
\forall v,w \in \c^n, \ \forall 1 \le k, l \leq n .
\ee
So by \eqref{eqc31}, there is an isometry $V$ that maps
each $p_{k,v}$ to $q_{k,v}$. If the span of the vectors $\{ p_{k,v} \}$ is
not all of $\c^{1 + \Np J}$, we can extend $V$ to the orthocomplement
so that it becomes an isometry (indeed, a unitary) from 
all of $\c^{1 + \Np J}$ to $\c^{1 + \Np J}$. 

With respect to the decomposition $\c \oplus \c^{J\Np}$,
write
\[
V \= 
\begin{bmatrix} A &  B\\
C & D\end{bmatrix}.
\]
We have
\be
\label{eqc32}
\begin{bmatrix} A &  B\\
C & D\end{bmatrix}
 [\pi_k \ot \idd_{\c^{1 + \Np J}} ] 
\begin{bmatrix}
\id{\cn} \\
[\d(\La)  \ot\id{{\c^{\Np}}} ]u(\La) 
\end{bmatrix}\=
[\pi_k \ot \idd_{\c^{1 + \Np J}} ] 
\begin{bmatrix}
p_0(\La) \\
u(\La) 
\end{bmatrix} .
\ee
Since \eqref{eqc32} holds for each $k$, we get that
\be
\label{eqc33}
\begin{bmatrix} \id{\cn} \ot A &  \id{\cn} \ot B\\
\id{\cn} \ot C & \id{\cn} \ot D\end{bmatrix}
\begin{bmatrix}
\id{\cn} \\
[\d(\La) \ot \id{{\c^{\Np}}} ]u(\La) 
\end{bmatrix}\=
\begin{bmatrix}
p_0(\La) \\
u(\La) 
\end{bmatrix} .
\ee

For $x$ in $\gdel \cap \mnd$,
define 
\beq
\lefteqn{\phi(x) \= \id{\cn} \ot A \ + \
 \left[(\idcn \otimes B) (\d(x) \ot \id{{\c^{\Np}}} )\right]
}
\\
&&
\left[ \idcn \otimes \id{{\c^{J\Np}}} - (\idcn \otimes D) (\d(x) \ot \id{{\c^{\Np}}} )\right]^{-1}
\idcn \otimes C.
\eeq
Then $\phi$ is in the unit ball of $H^\i(\gdel)$ by Theorem \ref{thmb2}. Moreover,
by\eqref{eqc33}, 
\[
\phi(\La) \= p_0(\La),
\]
as desired.
\ep

\section{The Nevanlinna Problem}
\label{secd}

There are two sources of non-uniqueness in the solution of the Pick interpolation problem.
 The first is the choice of $u$ in \eqref{eqc10};
the second is in the extension of $V$. This problem has been analyzed in the scalar case  by J.~Ball, T.~Trent and V.~Vinnikov \cite{btv01};  their ideas extend to our situation.

Let us suppose throughout this section that
\be
\label{eqd01}
\La \mapsto p_0(\La)
\ee
is a solvable Pick problem, and we have found a finite-dimensional space $\L$,
  an $\L(\c, \c^J \ot \L)$-valued nc polynomial $u$ satisfying 
\be
\label{eqd1}
 \idd_{\c^n} - p_0(\La)^* p_0(\La) \\
\=
u(\La)^* \left[ \idd_{\c^{nJ}\ot \L} - (\d(\La)^* \d(\La) \otimes \idd_{\L})
 \right] u(\La) ,
 \ee
and $V$ satisfying \eqref{eqc32}:
\be
\label{eqd2}
V  [\pi_k \ot \idd_{\c \oplus \c^J \ot \L} ] 
\begin{bmatrix}
\id{\cn} \\
[\d(\La) \ot \id{{\L}}   ]u(\La) 
\end{bmatrix}\=
[\pi_k \ot \idd_{\c \oplus \c^J \ot \L} ] 
\begin{bmatrix}
p_0(\La) \\
u(\La) 
\end{bmatrix} .
\ee
Let $\clo = \c \oplus \c^J \ot \L$, and 
\[
\cnt \ := \ 
\vee_{k=1}^n \vee_{v \in \cn} 
\begin{bmatrix}
\pi_k \, v \\
[(\pi_k \ot \id{{ \c^J \ot \L}})( \d(\La) \ot \id{{\c^{\Np}}} ) ]u(\La)  \, v
\end{bmatrix}  \subseteq \clo.
\]
Let 
\[
\cno
\ := \ 
\vee_{k=1}^n \vee_{v \in \cn} 
\begin{bmatrix}
\pi_k p_0(\La) \, v \\
(\pi_k \ot \id{{ \c^J \ot \L}}) u(\La) \, v
\end{bmatrix} \subseteq \clo ,
\]
and define $\cmt = \clo \ominus \cnt$ and $\cmo = \clo \ominus \cno$.
Then $V$ is an isometry from $\cnt$ onto $\cno$. Define a unitary 
\begin{eqnarray}
\nonumber
U : \cmo \oplus \cmt \oplus \cnt
 &\  \to \ &
\cmt \oplus \cmo \oplus \cno \\
\begin{bmatrix}
m_1 \\
m_2 \\n_2
\end{bmatrix}
& \mapsto &
\begin{bmatrix}
m_2 \\
m_1 \\ V n_2
\end{bmatrix}.
\label{eqd14}
\end{eqnarray}
By identifying $ \cmo \oplus \cmt \oplus \cnt \cong
\cmo \oplus \c \oplus  \c^J \ot \L$ and
$\cmt \oplus \cmo \oplus \cno
\cong
\cmt \oplus \c \oplus  \c^J \ot \L$, we can think of $U$ as a unitary from
$\c \oplus \cmo \oplus  \c^J \ot \L$ to $\c \oplus \cmt \oplus  \c^J \ot \L$, and it has a corresponding transfer function
$G$ that is a free $\L(\c \oplus \cmo, \c \oplus \cmt)$-valued  rational function
(since all the spaces are finite dimensional).
Write this $G$ as
\be
\label{eqd5}
 G \= \bordermatrix{&\c &\cmo \cr
\c &G_{11} & G_{12} \cr \cmt &G_{21}  & G_{22}} .
\ee
\bt
\label{thmd1} 
The function $\phi$ in $\bhigd$
satisfies $\phi(\La) = p_0(\La)$ 
 if and only if,
for some $u$ satisfying \eqref{eqd1} and $G$ the transfer function of $U$ in \eqref{eqd14}, there is a function
$\Theta$ in ${\rm ball}(H^\i_{\L(\cmo,\cmt)} (\gdel))$ such that
\be
\label{eqd6}
\phi \= G_{11} + G_{12} \Theta ( I_{\M_1} - G_{22} \Theta)^{-1}
G_{21} .
\ee
\et
\bp
($\Leftarrow$) This is a straightforward calculation.

($\Rightarrow$)
By Theorem~\ref{thmb2}, $\phi$ has a free $\d$-realization, and by Lemma~\ref{lemc2},
we can assume that $\{ u(x) : x \in \Vl \}$ lie in a finite dimensional space that
we can embed in $ \c^J \ot \L$. So we can assume that $\phi$ is the transfer function
of some unitary $X: \c \oplus \cj \ot \k \to \c \oplus \cj \ot \k$, and that $ \L \subseteq \k$.
For $x \in \gdel \cap \m^d_m$ we have
\be
\label{eqd8}
[ \idd_{\c^m} \ot X  ]
\begin{bmatrix}
 \idd_{\c^m} \\
(\d(x) \ot \idd_{\K}) \xi(x) 
\end{bmatrix} \, 
\=
\begin{bmatrix}
\phi(x)  \\
 \xi(x) 
\end{bmatrix} \,  .
\ee

Let $\kp = \k \ominus   \L$. Then 
\be
\label{eqd7}
 X \= \bordermatrix{&\cnt &\cmt \oplus \cj\ot \kp \cr
\cno &V & 0\cr \cmo \oplus \cj\ot\kp &0  & Y} .
\ee
Let $\Theta$ be the transfer function of $Y$. Then we claim that
\eqref{eqd6} holds.

Let $x \in \gdel \cap \m^d_m$ and $v \in \c^m$ be fixed for now.
Let
\begin{eqnarray}
\nonumber
p\=&  v \op 
( \d(x) \ot \id{\K} ) \xi(x) v  &\= 
n_2 \op m_2 \op h_2 \\
q \=&  
\phi(x) v \op  
  \xi(x) v  &\= 
n_1 \op m_1 \op h_1
\label{eqd11}
\end{eqnarray}
where $n_2 \in \cm\ot \cnt, m_2 \in \cm \ot \cmt, n_1 \in \cm \ot \cno,
 m_1 \in \cm \ot \cmo$ and $h_2,h_1 \in \cm \ot \cj \ot \kp$.
Note from \eqref{eqd8} that
\be
\label{eqd9}
[ \idd_{\cm} \ot X] p \= q .
\ee
Let $P^\prime$ be the projection from $\cj \ot \K$ to $\cj \ot \kp$. As $\d(x) \ot \id{\K}$ commutes
with $\id{{\c^m}} \ot P^\prime$, we get from \eqref{eqd11} that
\[
[\d(x) \ot \id{{\kp}}] h_1 \= h_2 .
\]
Therefore
\be
\label{eqd13}
[ \id{{\cm}} \ot Y ] ( m_2 \op [\d(x) \ot \id{{\kp}}] h_1 ) \= m_1 \op h_1 .
\ee
As $\Theta$ is the transfer function of $Y$, \eqref{eqd13} implies that
\be
\label{eqd15}
\Theta(x) \, m_2 \= m_1 .
\ee

Let $P$ be the projection from $\cm \ot \cj \ot \k$ onto $\cm \ot \cj \ot \L$, and
let $\eta = P \xi(x) v$.
Then under the identifications of $\cno \op \cmo$ 
and $\cnt \op \cmt$ with $ \c \op \cj \ot \L$, we get
\beq
n_1 \op m_1 &\=& \phi(x) v \op \eta \\
n_2 \op m_2 &\=& v \op (\d(x)  \ot \id{\L}) \eta .
\eeq
Then  from \eqref{eqd14}
\[
U : \ v \op m_1 \op (\d(x)  \ot \id{\L}) \eta \ \mapsto \
\phi(x) v \op m_2 \op \eta .
\]
By \eqref{eqd15} this gives
\be
\label{eqd16}
\begin{pmatrix}
G_{11} (x) & G_{12}(x) \\
G_{21}(x) & G_{22}(x) 
\end{pmatrix}
\begin{pmatrix}
v \\
\Theta(x) m_2 
\end{pmatrix}
\=
\begin{pmatrix}
\phi(x) v \\
 m_2 
\end{pmatrix}.
\ee
As \eqref{eqd16} holds for all choices of $x$ and $v$, we get \eqref{eqd6}, as desired.
\ep


\section{Extending functions defined on varieties}
\label{sece}

{\sc Proof of Theorem~\ref{thme1}:}
Without loss of generality, assume that
\be
\label{eqe5}
 \sup_{x \in \gv \cap \gdel} \| f(x) \| \= 1 .
\ee
Choose a sequence $(\l_j )_{j=1}^\i$ in $\gdel \cap \gv$ that is dense in the disjoint union topology,
so for all $\vare > 0$, for all $x \in \gdel \cap \gv$, there exists some $\l_j$ such that
$ \max_{1 \leq r \leq d} \| \l_j^r - x^r \| < \vare $.

Let $\La_n = \oplus_{j=1}^n \l_j$.
By Theorem~\ref{thmb1},  $f$ is 
 locally approximable by polynomials, and so has 
 the property that
\[
\forall x \in \gv \cap \gdel, \ f(x) \in \Mx.
\]
Therefore there is some polynomial $p_n \in \pd$ such that 
\be
\label{eqe3}
p_n(\La_n) \= f(\La_n) .
\ee
Moreover, if $x \in  \gv \cap \gdel$, then by Theorem~\ref{thmb1} again, one can approximate $f$
at $x \oplus \La_n$ by a sequence of free polynomials, and so by Lemma~\ref{lemb1}
\be
\label{eqe4}
\forall x \in \gv \cap \gdel, \ f(x) = p_n(x).
\ee
As $\Vl \subseteq \gv$, putting \eqref{eqe3}, \eqref{eqe4} and \eqref{eqe5} together, the hypotheses of Theorem~\ref{thmpick} are satisfied, so there exists $\phi_n$ in $\hinfo$ such that
\[
\phi_n(\La_n) \= f(\La_n) .
\] 
By Theorem~\ref{thmb3}, some subsequence of $\phi_n$ converges locally uniformly (in the disjoint union topology)
 to a function
$\phi$ in $\hinfo$. Moreover, for each $j$, $\phi(\l_j) = f(\l_j)$, so by continuity, $\phi$ is an extension of $f$.
\ep

\section{Examples}
\label{secf}

\addtocounter{equation}{1}
{\bf Example \theequation}
Let $q_1, \dots , q_m$ be polynomials in $d $ commuting variables, and let $V = \{ z \in \c^d : q_i(z) = 0 , i = 1,\dots,m \} .$
Let $f$ be a (scalar-valued)  holomorphic function defined on $ V \cap \D^d$.

Let $T$ be a $d$-tuple of commuting matrices that are strict contractions, 
and such that $q_i(T) = 0 $ for $i =1,\dots,m$.
If they are simultaneously diagonizable,
then their joint eigenvalues lie in $V \cap \D^d$, and  it makes sense to define $f(T)$ by
applying $f$ to the diagonal entries, in the basis of joint eigenvectors.
If the matrices are not simultaneously diagonizable, then one can still define $f(T)$, 
either by the Taylor functional calculus \cite{tay70b}, or, more constructively, as in
\cite{am13}.

Let us write $\cal F$ for the set of all $T = (T^1, \dots , T^d)$ of commuting matrices such that
$q_i(T) = 0, i = 1,\dots,m$, and such that $\| T \| < 1$. 
Note that ${\cal F} =  \gv \cap \gdel$, where $\ga$ is the algebra generated by
$q_1, \dots, q_m$ and the polynomials $\{ x^i x^j - x^j  x^i: 1 \leq i < j \leq d \}$,
$\gv = \varga$,
and $\d(x)$ is the diagonal matrix with entries $x^1, x^2, \dots, x^d$.
Define a norm on holomorphic functions on 
$ V \cap \D^d$
by
\[
\| f \|_{\gv \cap \gdel} \= 
\sup\{ \| f(T) \| : T \in {\cal F} \}.
\]
To apply Theorem~\ref{thme1}, we need to know that $f$ is a free holomorphic function on
$\gv \cap \gdel$, in other words that locally in $\gdel$ it extends to a free holomorphic function
(\ie it can be applied to non-commuting matrices). This is true, and is proved in \cite{am13}.
Then Theorem~\ref{thme1}
 asserts that there is a bounded extension $\phi$ of $f$, 
defined on the set $\{ R \in \m^d : \| R \| < 1 \}$, if and only if 
$\| f \|_{\gv \cap \gdel}$ is finite. Moreover, if this quantity is finite, then $\phi$ can be found with exactly this norm.
In particular, an extension to the non-commuting ball $\gdel$ can always be found with the same norm as is attained
by evaluating on commuting matrices in the variety.
\vs
\addtocounter{equation}{1}
{\bf Example \theequation}
Specializing the previous example to the case $d=2$,
and using 
And\^o's inequality \cite{and63}, we conclude the following: if we wish to extend a polynomial
$p_0$ off $V \cap \D^2$,  where $V$ is the joint sero set of the $q_i$'s, then the minimum norm of the extension $\phi$ is the same when calculated
as a scalar-valued function in $H^\i(\D^2)$, as a function on pairs of commuting contractive matrices,
or as a function on pairs of contractive matrices. The norm is attained, and is given by
\begin{eqnarray}
\nonumber
\sup_{n \in \n}\  \sup \{ \| p_0(T) \| \ : \
T \in \mn^2, \, \|T^1 \| < 1, \| T^2 \| < 1, \\
 \, T^1 T^2 = T^2 T^1, \ q_i(T) = 0 \, \forall \, 1 \leq i \leq m \} .
\label{eqfx1}
\end{eqnarray}
Unless $V \cap \D^2$ is a retract of $\D^2$, one can  by
\cite{agmc_vn} always find some $p_0$ so that 
\eqref{eqfx1} is strictly greater than 
$$\sup \{ | p_0(z) | \ : \ z \in \D^2 \cap V \} .$$
\vs
\addtocounter{equation}{1}
{\bf Example \theequation}
Suppose $\delta(x)$ has first column $x^1, \dots, x^d$  and its other entries zero, so $$
\gdel \= 
\{ T \ : \ T^{1*} T^1 + \cdots T^{d*} T^d < 1 \} .
$$
(This is called the row ball).
Suppose $\La, H \in \gdel \cap \mnd$ and one wishes to solve the interpolation problem
\begin{eqnarray}
\nonumber
\phi(\La) &\= & W \\
\label{eqf1}
D \phi(\La) [H] &\=& X ,
\end{eqnarray}
where $D \phi(\La) [H]$, the derivative of $\phi$ at $\La$ in the direction $H$, is defined by
\[
D \phi(\La) [H]
\=
\lim_{t \to 0} \frac{\phi(\La + t H) - \phi(\La)}{t} .
\]
A necessary condition to find a function $\phi \in \hinfgdel$ solving this problem is
that there is some free polynomial $p_0$ with $p_0(\La) = W$ and
$D p_0(\La) [H] = X$. 
The minimum norm of a solution can be found from Theorem~\ref{thme1} by letting
\[
\ga \= 
\{ p \in \pd \ : \
p(\La) = 0,\ Dp (\La) [H] = 0 \} ,
\]
$\gv = \varga$, and calculating 
\[
\sup_{x \in \gv \cap \gdel} \| p_0(x) \|  .
\]
The problem can also be solved using Theorem \ref{thmpick}, as 
\eqref{eqf1} is the same as solving the one point problem
\[
\begin{pmatrix}
\La & H \\
0 & \La
\end{pmatrix}
\ \mapsto \
\begin{pmatrix}
W & X \\
0 & W
\end{pmatrix},
\]
since by \cite[Prop 2.5]{hkm11b}, for any continuous nc-function $f$, one has
\[
f
\begin{pmatrix}
\La & H \\
0 & \La
\end{pmatrix}
\=
\begin{pmatrix}
f(\La) & Df(\La)[H] \\
0 & f(\La)
\end{pmatrix}.
\]

\bibliography{../../references}

\end{document}